\documentstyle{article}
 \def\vt{t\kern-0.22em\raise.18ex\hbox{\char'47}\lower.18ex\hbox{}\kern-0.08em}
 
\newtheorem{th}{Theorem}[section]

\newtheorem{con}{Conjecture}[section]

\newcommand{\old}[1]{{}} 
\newcounter{obr}
\newcounter{tabul}
\begin{document}
\title{Strong Pseudo Transitivity and Intersection Graphs
}
\author{Farhad Shahrokhi\\
Department of Computer Science and Engineering,  UNT\\
farhad@cs.unt.edu
}

\date{}
\maketitle
\thispagestyle{empty}
\date{} \maketitle
\begin{abstract} 
A directed graph $G=(V,E)$ is {\it strongly pseudo transitive}  if there is a partition $\{A,E-A\}$ of $E$  so that graphs $G_1=(V,A)$ and 
$G_2=(V,E-A)$ are transitive,  and additionally, if  $ab\in A$ and $bc\in E $ implies that  $ac\in E$.  A strongly pseudo  transitive graph $G=(V,E)$ is strongly pseudo transitive of the first type, if $ab\in A$ and $bc\in E$ implies $ac\in A$.  
An undirected graph is co-strongly pseudo transitive (co-strongly pseudo transitive of the first type)  if its  complement has an orientation which is  strongly pseudo transitive (co-strongly pseudo transitive of the first type). Our purpose is  show that  the  results  in  computational geometry \cite{CFP, Lu} and intersection graph theory  \cite{Ga2, ES}  can be unified and extended,  using  the  notion of strong pseudo transitivity.  As a  consequence  the   general algorithmic  framework in  \cite{Sh} is applicable to solve the maximum independent set in $O(n^3)$ time in a variety of problems, thereby, avoiding case by case lengthily arguments for  each problem.

We show that  the  intersection graphs of axis parallel rectangles intersecting a diagonal line from bottom, and half segments  are co-strongly pseudo transitive.
In addition,  we show that  the class of the interval filament graphs is co-strongly transitive  of the  first type, and hence  the class of  polygon circle  graphs which is contained in the class of interval filament graphs (but contains the   classes of chordal graphs, circular arc, circle, and outer planar graphs), and the class of incomparability graphs are strongly transitive of the first type.  For class of chordal graphs we give two different proofs, using two different characterizations, verifying that they are  co-strongly transitive  of the  first type. Furthermore, we note that the class of co-strongly pseudo transitive  graphs of the first type  is properly contained in the class co-pseudo  transitive graphs, unless $P=NP$, and the class of tree filament graphs is  properly contained in the class co-pseudo  transitive graphs, unless $P=NP$. 
Computational consequences are presented. 
 \end{abstract}

\section{Introduction}
``The first genuine monograph on graph theory (K\"onig, 1936)
had the following subtitle: {\it Combinatorial Topology of
Systems of Segments} \cite{K}. Although graph theory and topology
stem from the same root, the connection between them has
somewhat faded away in the past few decades. In the most
prolific new areas of graph theory including Ramsey theory, extremal
graph theory and  random graphs, graphs are regarded as
abstract binary relations rather than systems of segments.
It is quite remarkable that {\em traditional} graph theory is
often incapable of providing satisfactory answers for the
most natural questions   concerning the  drawings of graphs'' \cite{Pa}.

Recall that a directed graph $G=(V,E)$ is {\it transitive} if $ab\in E$ and $bc\in E$ implies $ac\in E$. One can view a transitive graph as an alternate way of defining a partial order 
\cite{Tr2}. By dropping the orientation on the edges of a transitive  graph we obtain a comparability graph. The complement of a comparability graph is an incomparability
graph. In \cite{Sh} we  introduced the  concepts of {\it pseudo  transitivity} and {\it strong pseudo  transitivity}  which are  ways of generalizing the concept of transitivity in graphs. The motivation behind writing this  paper arose from recent results in computational geometry \cite{CFP, Lu} and  intersection graphs theory \cite{ES}, \cite{Ga2}  that can be unified  and generalized  using the  concept of strong pseudo  transitivity.

Let $S$ be a finite set and $\{S_1,S_2,...,S_k\}$ be a collection of subsets of $S$.
The intersection (overlap) graph of $S$ is a graph with  the vertex set $S$, where $S_iS_j$ is an edge if $S_i\cap S_j\ne \emptyset$ 
($S_i\cap S_j\ne \emptyset$, and $S_i\not\subseteq S_j$,  $S_j\not\subseteq S_i$). 
  
Many intersecting graphs are the intersection or the overlap graphs of combinatorial or geometric structures. For instance  incomparability graphs are intersection graphs of $x-$monotone curves intersecting two lines parallel  to the  $y$ axis \cite{GR},  overlap graphs are  overlap graphs of intervals on a line, chordal graphs are the intersection graph of subtrees of a tree \cite{Ga1},  polygon-circle graphs are  the intersection graphs of  convex polygons in a circle, Interval  filament  graphs are the intersection graphs of filaments over intervals \cite{Ga2} and,
subtree filament are intersection graphs of filaments over trees \cite{Ga2}, or overlap graph of subtrees of a tree \cite{ES}.

 A directed graph $G=(V,E)$ is {\it  pseudo transitive} if there is a partition $\{A,E-A\}$ of $E$  so that graphs $G_1=(V,A)$  is transitive,  and additionally, if  $ab\in A$ and $bc\in E $ implies that  $ac\in E$.
 A directed graph $G=(V,E)$ is {\it strongly pseudo  transitive} if there is a partition $\{A,E-A\}$ of $E$  so that graphs $G_1=(V,A)$ and 
$G_2=(V,E-A)$ are transitive,  and additionally, if  $ab\in A$ and $bc\in E $ implies that  $ac\in E$. 
Let   $G=(V,E)$  be pseudo  transitive (pseudo transitive) with the underlying partition $\{A, E-A\}$, then $G$ is  {\it  pseudo transitive  of the first type} (pseudo transitive of the first type)  if  $ab\in A$ and $bc\in E$ implies that  $ac\in A$.
An undirected graph $G$ is co-pseudo transitive (co-pseudo transitive of the first type)  if the complement of $G$  has an orientation that is pseudo transitive (pseudo transitive of the first type one. $G$ is  co-strongly pseudo transitive (co-strongly transitive  of the first type), if  the complement of $G$    has an  orientation which  is strongly pseudo transitive (co-strongly transitive of the first type). Co-pseudo transitive graphs contain the intersection graph of many geometric structures.
For instance,  the following result was shown  in \cite {Sh}, with  a slightly altered  language. 
\begin {th}\label{t0}
{\sl 
Let $P$ be a finite collection of bounded closed subsets of $R^k$, then 
the intersection graph of $P$ is co-pseudo transitive of the first type. 
}

\end{th}

A half segment is a straight  line segment that has  one end point on the $x$-axis,  another end point in the upper half plane, and makes  an  acute angel with $x-$axis. 
Motivation behind introducing these segments  arose from the work of Pach and Torocsik \cite{PT}  on  geometric  graphs. Biro and Trotter \cite{BT} studied  properties  of partial orders arising from half segments. Computing the maximum independent set in the intersection graph of a set of rectangles is a fundamental problem arising  in map labeling \cite{TDE}. 
 Since the general  version of this problem is known to be NP-hard, some researchers have  focused to solve the special version of the problem, including the instances where all rectangles are intersected by a diagonal line from below.
See the work of Lubiw \cite{Lu}, and  Correa, Feuilloley,Perez-Lantero , and  Soto \cite{CFP}, which provide an $O(n^3)$ and $O(n^2)$  time algorithm, respectively.

In this paper we explore the connections between the  graph  classes mentioned above  and  the class of  co-strongly pseudo transitive graphs.  Specifically we show  that  the  intersection graphs of  half line segments and  axis parallel rectangles intersecting a diagonal line from bottom are co-strongly pseudo transitive. Moreover,  we show that  the class of the interval filament graphs is co-strongly transitive of the first type, and hence  the class of  polygon circle  graphs which is contained in the class of interval filament graphs (but contains the   classes of chordal, circular arc, circle, and outer planar graphs), and the class of incomparability graphs are strongly transitive of the first type.  For the  class of chordal graphs which is  contained in the class of polygon circle graphs, we provide two  direct direct proofs,  showing  that they are co-strongly transitive of the first type. Additionally, we present some results concerning the 
Containment of different classes.  A contribution of our work is  to connect and unify the problems  in  computational geometry \cite{CFP, Lu},  intersection graph theory  \cite{Ga2, ES} and combinatorics 
 \cite{BT}  using  the  notion of strong pseudo transitivity, thereby, showing they are all   amenable to the  algorithmic  framework in  \cite{Sh}  for  solving  the maximum independent set in $O(n^3)$ time, thereby, avoiding case by case or lengthy arguments for each scenario.

\section{Structural Results}
\begin{th}
{\sl Let $R$ be a set of axis parallel  rectangles in the plane all of them are intersected by a diagonal  line $l$ with the property that if two elements of $R$ intersect, then they also intersect below $l$.  Let  $G=(V,E)$ be the intersection graph of $S$,   then $G$ is co-strongly pseudo transitive.}
 \end{th}
 {\bf Proof.} For any $C\in R$, let $x_C$ and $y_C$ denote the smallest  $x$ coordinate, and smallest  $y$ coordinate of four corners of $C$, respectively. Now let $W, Z\in R$ with $WZ\notin E$,  so   that  $x_W\le x_Z$
 (the case $x_W>x_Z$ is symmetric). If $y_W\le y_Z$, then, orient $WZ$ from $W$ to $Z$, and place it in $A$. Otherwise if $y_W>y_Z$, then, still orient $WZ$ from $W$ to $Z$ but place it in $B$. It can be verified   that $A\cap B=\emptyset $, and that any non edge of $G$ has a orientation  in $\hat E=A\cup B$. 
 Furthermore, it can be shown that  this case, the directed acyclic  graph   $H=(V,\hat E)$   is strongly pseudo-transitive, with the partition $\{A, \hat E-A\}$. $\Box$.
 
The following result was mentioned in \cite{Sh} without an explicit proof. Next,  we specifically state and prove it. 
 
 \begin{th}
 {\sl Let $R$ be a set of half segments in the plane, and let  $G=(V,E), V=R$ be the intersection graph of $R$,   then $G$ is co-strongly pseudo transitive.}
 \end{th}
 {\bf Proof.} For any $r\in R$, let $x_1(r)$ and $x_2(r)$ denote the smallest $x$ coordinate, and the largest $x$ coordinate of any points in $r$. 
 Now let $r,s\in R$ so that $rs\notin E$ so that $x_1(r) < x_1(s)$.  if 
 $x_2(s)\ge x_2(r)$, then orient $rs$ from $r$ to $s$ and place $rs$ in $A$.
 Otherwise, if $x_2(s)<x_2(r)$ (note the assumption $rs\notin E$), then still  orient $rs$ from $r$ to $s$, but  place the directed  edge $rs$ in $B$. The remaining of the proof copies previous the theorem. $\Box$ 
 
 \begin{th}
{\sl Let $G=(V,E)$ be an interval filament graph, then $G$ is co-strongly pseudo transitive of the first type.}
 \end{th}

 {\bf Proof.} Consider a representation of $G$, where  $I$ is  a set of intervals on the real line, and for each $i\in I, C_i$ is a collection of curves on the half plane above $i$ that connects the end points of $i$. Then $V=\cup_{i\in I}C_i$ and furthermore $xy\in E$, if  $x,y\in C_i$ for some $i\in I$, and  $x$ and $y$  intersect.  
 Now let  $x,y\in V, xy\not\in E$.  If $x\in C_i, y\in C_j, i<j$, then orient $xy$ from $x$ to $y$ and place $xy$ in $A$. Otherwise $x,y\in C_i$ for some $i\in I$. Since  $x$ and $y$ do not intersect and connect the endpoint for interval $i$, the area under one of them (say $x$) contains the area under the other (say $y$). In this case orient $xy$ from $x$ to $y$ and place it in $B$. Note that  $A\cap B=\emptyset$, every $xy\notin  E$ has an orientation in $\hat E=A\cap B$. Furthermore the  directed acyclic  graph   $H=(V,\hat E)$, is strongly pseudo transitive, and  for any  $xy\in A$ and $yz\in  \hat E$, we have  $xy\in A$. $\Box$.

\begin{con}\label{c1}
 {\sl The class of interval filament graphs is properly contained in the class co-strongly pseudo  transitive of the first type.}
 
  \end{con}

Since chordal graphs play an important role in graph theory,  we 
give  two different direct proofs, based on  different representations,   showing that  they are co-strongly pseudo transitive. The first proof  uses the  characterization that every chordal graph is the intersection graph of subtrees of the a tree, where, the second assumes  a perfect elimination ordering is given.

 Let $L={v_1,v_2,...,v_n}$ be a perfect elimination ordering  (PEO) of a chordal graph $G=(V,E)$. A {\it canonical depth first search tree}  of $G$ (with respect to $L$) is  a depth first search  spanning tree rooted at $v_1$  constructed applying  the following simple rule for visiting vertices: Assume vertex  $v_i$ is  currently visited, then,  select   the next vertex to visit  to be the smallest indexed vertex $v_j, j>i$ among all unvisited vertices adjacent to $v_i$ in $G$.

   \begin{th}\label{t1}
{\sl Every choral graph $G=(V,E)$ is co-strongly pseudo transitive of the first type.}
 \end{th}
 {\bf First Proof.} 
  Let $T$ be a tree,  let $V=\{T_1,T_2,...,T_k\}$ be a set of subtrees of $T$. Assign a root $r$ to $T$, embed $T$ in the  plane, and assign  root $r_i$ to each $T_i\in V$, which is the closest vertex of $T_i$ to $r$. Let  $ G=(V,E)$ be the intersection graph of these subtrees, let ${\bar G}=(V,{\bar E})$ be the complement of $G$.  To prove the claim we will show there is a suitable  orientation on  $\bar E$.  
  For any  $e=T_iT_j\in\bar E$ so that $r_i$ is not an ancestor of $r_j$ in $T$, and $r_i$  to the left of $r_j$ in the planar  embedding of $T$, orient $e$ from $T_i$ to $T_j$ and place the resulting oriented  edge in $A$.  For any $e=T_iT_j\in {\bar E}$ so that  $r_i$ is an ancestor  of $r_j$ in $T$
 orient  $e$ from $T_i$ to $T_j$ and place the resulting  oriented  edge in $B$.
 
 It can be verified that  $A\cap B=\emptyset$, every $e\in {\bar E}$ has an orientation in $\hat E=A\cap B$, and  that directed acyclic  graph   $H=(V,\hat E)$,   is strongly pseudo  transitive, thus verifying the claim.  Moreover, in this case $xy\in A$ and $yz\in  \hat E$ implies $xy\in A$.

{\bf Second Proof.} 
Let $L=\{v_1,v_2,....,v_n\}$ be a PEO  of $G$. Let $T$ be a canonical depth first search spanning tree of $G$ rooted at $v_1$. Let ${\bar e}=xy\not\in E$ with $dfs(x)<dfs(y)$. If $x$ and $y$ are on two different branches of $T$ then orient ${\bar e}$ from $x$ to $y$ and place $xy$ in $A$. Otherwise,  $x=v_i$ and $y=v_j$ are on the same  branches of $T$. Observe  in this case that  we must have $i<j$.   Since $T$ is a canonical depth first tree, and, orient ${\bar e}$ from $x$ to $y$ and place $xy$ in $B$.

{\bf  Claim.} Let $v_i,v_j, v_k, i<j<k$ be three  vertices on the same branch of $T$. If $v_iv_j\notin E$, and $v_jv_k\notin E$, then $v_iv_k\notin E$.

It can verified (using the claim and properties of $T$) that $A\cap B=\emptyset$, every $e\in {\bar E}$ has an orientation in $\hat E=A\cap B$,  and that  directed acyclic  graph   $H=(V,\hat E)$   is strongly pseudo  transitive. Additionally,  for any  $xy\in A$ and $yz\in  \hat E$, we have  $xy\in A$.
$\Box$

We finish this section by the  establishing  some containment properties.

\begin{th}\label{t10}
{\sl
 
$(i)~$The class of co-pseudo  (co-strongly pseudo) transitive graphs of the first type is contained in the class of co-pseudo  (co-strongly pseudo) transitive graphs.

$(ii)~$The class of co-strongly pseudo transitive  graphs of the first type  is properly contained in the class co-pseudo  transitive graphs of the first type, unless $P=NP$.

$(iii)~$The class of tree filament graphs is  contained in the class 
co-pseudo  transitive graphs.

$(iv)~$The class of tree filament graphs is  properly contained in the class 
co-pseudo  transitive graphs, unless $P=NP$.
}
\end{th}

{\bf Proof.} Clearly $(i)$ holds. For $(ii)$, first note that by Theorem \ref{t0}
the intersection graph  $G$ of a set of rectangles is co-pseudo transitive of the first time. Next note that  computing the maximum independent set of $G$ is NP-hard, but can be done in $O(n^3)$ time for any co-strongly pseudo transitive graph. We omit proof of $(iii)$. $(iv)$ follows that graphs of boxicity two are co-pseudo transitive and  computing their maximum independent set is known to be NP-hard, but  computing 
maximum independent set in tree filament graphs can be done in polynomial time.
$\Box$.

\section{Algorithmic Consequences} 
The following result was shown in \cite{Sh}

\begin{th}
{\sl 
Let $H=(V,F)$ be strongly pseudo-transitive. The maximum weighted chain  can be computed in $O(\sum_{x\in V}deg^2(x)+n^2)$. 
}
 \end{th} 
 Using our notations the above theorem implies.
 \begin{th}
{\sl 
Let $G=(V,E)$ be co-strongly pseudo-transitive. The maximum weighted independent set  can be computed in $O(n^3)$. 
}
 \end{th} 

The above theorem  implies the following general result. 
\begin{th}
{\sl 
Let $G=(V,E)$ be one of the following graphs $(i)$ incomparability  $(ii)$ overlap , $(iii)$ chordal, $(iii)$ polygon circle $(iv)$ interval filament $(v)$, or, 
$(vi)$ intersection graph half segments,  $(vii)$ intersection  graph axis parallel rectangles intersected by a diagonal line. Then,   the maximum weighted independent set  can be computed in $O(n^3)$. 
}
 \end{th}

Note that  our  general frame work  extends the work of  Lubiw 
 \cite{Lu} who showed the weighted maximum independent set of rectangles all which  have their right most corner on a line  can be computed in  $O(n^3)$ time, but gives  weaker result than a  more recent work of  Correa, Feuilloley,Perez-Lantero ,  Soto 
 \cite{CFP} that had  $O(n^2)$  time complexity.





\begin{thebibliography}{99}

\bibitem{BT}
Biro, C., Trotter, W. T., Segment Orders, DCG, vol. 43, no. 3, pp. 680-704, 2010.

\bibitem{CFP}
Correa, Feuilloley, Perez-Lantero,  Soto, Independent and Hitting Sets 
of Rectangles Intersecting a Diagonal Line: 
Algorithms and Complexity, Discrete and  Computational Geometry, 
volume 53, pages 344-365 (2015).  

\bibitem{Go} Golumbic  M. C., Algorithmic Graph Theory and Perfect Graphs (Annals of Discrete Mathematics, Vol 57, North-Holland Publishing Co., Amsterdam, The Netherlands, 2004.

\bibitem{GR}
Golumbic, M.; Rotem, D.; Urrutia, J. (1983), "Comparability graphs and intersection graphs", Discrete Mathematics, 43 (1): 37–46.
 
\bibitem{K} 
K\"onig D.: {\it Theorie der endlichen und
unendlichen Graphen,} Leipzig, 1936. English translation:
{\it Theory of Finite and Infinite Graphs}, Birkh\"auser
Verlag, Boston, 1990.

\bibitem{ES}
 Enright J.,  Stewart L., 
Equivalence of the filament and overlap graphs of subtrees of limited trees,
	arXiv:1508.01441 [cs.DM], 2017. 
	
\bibitem{Lu} 
Lubiw, A.,  A weighted min-max relation for intervals. J. Comb. Theory, Ser. B 53
(1991) 151–172

\bibitem{Ga2}
 Gavril, F.,
Maximum weight independent sets and cliques in intersection graphs of filaments.
Information Processing Letters 73(5-6) 181-188 (2000).

\bibitem{Ga1}
Gavril,  F.,
The intersection graphs of subtrees in trees are exactly the chordal graphs, Journal of Combinatorial Theory, Series B, Volume 16, Issue 1, February 1974, Pages 47-56.

\bibitem{Pa}
 Pach, J. CBMS-NSF  conference, UNT, Denton Tx, 2002.

\bibitem{Sh}
 Shahrokhi, F. Algorithms For Longest Chains In Pseudo- Transitive Graphs, 
Congressus Numerantium, 221 (2014), p 21-30,  arXiv:1701.05286v1 [cs.CG].

\bibitem{Tr} Trotter W.T.,   New perspectives on interval orders and interval
graphs, in Surveys in Combinatorics, Cambridge Univ. Press, 1977, 237-286.

\bibitem{Tr2} 
Trotter,  W. T.,
Combinatorics and Partially Ordered Sets: Dimension Theory (Johns Hopkins Studies in the Mathematical Sciences) December 18, 2001.



\bibitem{PT}
Pach J., and Torocsik, J., Some geometric applications of Dilworth’s theorem, Disc. Comput. Geometry 21(1994). 83-95.

\bibitem{TDE}
Tollis, G.,  Di Battista, G.,  Eades, P. ,  Tamassia,  R.,
Graph Drawing: Algorithms for the Visualization of Graphs 1st Edition
by Ioannis  
\end{thebibliography}
\end{document}